\documentclass[12pt]{article}
\usepackage{amsxtra,amssymb,amsthm,amsmath,latexsym}

\theoremstyle{plain}
\newtheorem{definition}{Definition}[section]
\newtheorem{theorem}{Theorem}[section]

\newtheorem{corollary}{Corollary}[section]
\newtheorem{remark}{Remark}[section]

\newtheorem*{remark3.2}{Remark 3.2}
\newtheorem*{theorem3.3}{Theorem 3.3}

\def\nd{\noindent}

\def\R{{\mathbb R}}

\def\C{{\mathbb C}}
\def\oH{\buildrel\circ\over H}
\def\oH1{\buildrel\circ\over H\kern-.02in{}^1}
\def\qed{{\hfill $\Box$}}

\def\supp{\hbox{\,supp\,}}
\def\const{\hbox{\,const\,}}

\begin{document}


\title{An essay on some problems of approximation theory
   \thanks{key words: stable differentiation, approximation theory,
   property $C$, elliptic equations, Runge-type theorems, scattering solutions
    }
   \thanks{Math subject classification: 41-XX, 30D20, 35R25, 35J10, 35J05,
   65M30 }
}

\author{
A.G. Ramm\\
 Mathematics Department, Kansas State University, \\
 Manhattan, KS 66506-2602, USA\\
ramm@math.ksu.edu\\
}

\date{}

\maketitle\thispagestyle{empty}

\begin{abstract}
Several questions of approximation theory are discussed: 1) can one
approximate stably in $L^\infty$ norm $f^\prime$ given approximation
$f_\delta, \parallel f_\delta - f \parallel_{L^\infty} < \delta$,
of an unknown smooth function $f(x)$, such that
$\parallel f^\prime (x) \parallel_{L^\infty} \leq m_1$?

2) can one approximate an arbitrary $f \in L^2(D), D \subset \R^n, n \geq 3$,
is a bounded domain, by linear combinations of the products
$u_1 u_2$, where $u_m \in N(L_m), m=1,2,$ $L_m$ is a formal linear partial
differential operator and $N(L_m)$ is the null-space of $L_m$ in $D$,
$N(L_m) := \{w: L_m w=0 \hbox{\ in\ } D\}$?

3) can one approximate an arbitrary $L^2(D)$ function by an entire function
of exponential type whose Fourier transform has support in an arbitrary
small open set? Is there an analytic formula for such an approximation?
\end{abstract}


\section{Introduction}
In this essay I describe several problems of approximation theory which I
have studied and which are of interest both because of their mathematical
significance and because of their importance in applications.

\subsection{} 
The first question I have posed around 1966. The question is: suppose
that $f(x)$ is a smooth function, say $f \in C^\infty (\R)$, which
is $T$-periodic (just to avoid a discussion of its behavior near the
boundary of an interval), and which is not known; assume that its
$\delta$-approximation $f_\delta \in L^\infty (\R)$ is known,
$\| f_\delta - f \|_\infty < \delta$, where $\| \cdot \|_\infty$ is the
$L^\infty (\R)$ norm. Assume also that
$\| f^\prime \|_\infty \leq m_1 < \infty$. Can one approximate stably in
$L^\infty (\R)$ the derivative $f^\prime$, given the above data
$\{\delta, f_\delta, m_1\}$?

By a possibility of a stable approximation (estimation) I mean the existence
of an operator $L_\delta$, linear or nonlinear, such that
$$\sup_{\substack{f \in C^\infty (\R) \\
        \|f - f_\delta \|_\infty \leq \delta, \| f^\prime \|_\infty \leq m_1}}
        \|L_\delta f_\delta - f^\prime \|_\infty \leq \eta (\delta) \to 0
        \hbox{\ as\ } \delta \to 0,
        \eqno{(1.1)} $$
where $\eta(\delta) >0$ is some continuous function, $\eta (0) = 0$.
Without loss of generality one may assume that $\eta(\delta)$ is monotonically
growing.

In 1962-1966 there was growing interest to ill-posed problems. Variational
regularization was introduced by D.L. Phillips \cite{P} in 1962 and a year
later by A.N. Tikhonov \cite{T1} in 1963. It was applied in \cite{DI} in
1966 to the problem of stable numerical differentiation. The method for
stable differentiation proposed in \cite{DI} was complicated.

I then proposed and published in 1968 \cite{R1} the idea to use a divided
difference for stable differentiation and to use the stepsize $h=h(\delta)$
as a regularization parameter. If
$\| f^{\prime \prime} \| \leq m_2$, then
$h(\delta) = \sqrt{\frac{2\delta}{m_2}}$, and if one defines (\cite{R1}):
$$L_\delta f_\delta :=
        \frac{f_\delta (x+h(\delta)) - f_\delta (x-h(\delta))}{2h(\delta)},
        \quad h(\delta) = \sqrt{\frac{2\delta}{m_2}},
        \eqno{(1.2)} $$
then
$$\| L_\delta f_\delta - f^\prime (x) \|_\infty \leq
        \sqrt{2m_2 \delta} := \varepsilon (\delta).
        \eqno{(1.3)}$$

It turns out that the choice of $L_\delta$, made in \cite{R1}, that is,
$L_\delta$ defined in (1.2), is the best possible among all linear and
nonlinear operators $T$ which approximate $f^\prime (x)$ given the
information $\{\delta, m_2, f_\delta\}$. Namely, if
${\cal K} (\delta, m_j) := \{f:f \in C^j (\R), m_j < \infty,
  \| f-f_\delta\|_\infty \leq \delta\}$,
and $m_j = \| f^{(j)} \|_\infty$, then
$$\inf_T \sup_{f \in {\cal K}(\delta, m_2)}
  \| Tf_\delta - f^\prime \| \geq \sqrt{2m_2 \delta}.
  \eqno{(1.4)} $$
One can find a proof of this result and more general ones in \cite{R2},
\cite{R3}, \cite{R6}, \cite{R7} and various applications of these results
in \cite{R2} - \cite{R5}.

The idea of using the stepsize $h$ as a regularization parameter became
quite popular after the publication of \cite{R1} and was used by many
authors later.

In \cite{R25} formulas are given for a simultaneous approximation of
$f$ and $f^\prime$.

\subsection{} 
The second question, that I will discuss, is the following one: can
one approximate, with an arbitrary accuracy, an arbitrary function
$f(x) \in L^2(D)$, or in $ L^p(D)$ with $p\geq 1$, by a linear
combination of the products $u_1u_2$,  
 where $u_m \in N(L_m)$, $m=1,2,$ $L_m$ is a formal linear partial
differential operator, and $N(L_m)$ is the null-space of $L_m$ in $D$,
$N(L_m) := \{w: L_m w=0 \hbox{\ in\ } D\}$? 

This question has led me to the notion of
property $C$ for a pair of linear formal partial differential operators
$\{L_1, L_2\}$. 

Let us introduce some notations.
Let $D \subset R^n, n \geq 3$, be a bounded domain,
$L_m u(x) := \sum_{|j| \leq J_m} a_{jm} (x) D^j u(x), m = 1,2,$ $j$ is a
multiindex, $J_m \geq 1$ is an integer, $a_{jm} (x)$ are some functions
whose smoothness properties we do not specify at the moment,
$D^j u = \frac{\partial^j u}{\partial x_1^{j_1} \dots \partial x^{jn}_n}$,
$|j| = j_1 + \dots + j_n$. Define
$N_m := N_D (L_m) := \{w : L_m w=0 \hbox{\ in\ } D\}$,
where the equation is understood in the sense of distribution theory. Consider
the set of products $\{w_1 w_2\}$, where $w_m \in N_m$ and we use all the
products which are well-defined. If $L_m$ are elliptic operators and
$a_{jm} (x) \in C^\gamma (\R^n)$, then by elliptic regularity the functions
$w_m \in C^{\gamma + J_m}$ and therefore the products $w_1 w_2$ are well
defined.

\begin{definition} 
A pair $\{L_1, L_2\}$ has property $C$ if and only if the set
$\{w_1 w_2\}_{\forall w_m \in N_m}$ is total in
$L^p (D)$ for some $p \geq 1$.

In other words, if $f \in L^p (D)$, then
$$\left\{ \int_D f(x) w_1 w_2 dx = 0, \quad \forall w_m \in N_m \right\}
  \Rightarrow f=0, \eqno{(1.5)} $$
where $\forall w_m \in N_m$ means for all $w_m$ for which the products
$w_1 w_2$ are well defined.
\end{definition}

\begin{definition} 
If the pair $\{L,L\}$ has property $C$ then we say that the operator $L$ has
this property.
\end{definition}
From the point of view of approximation theory property $C$ means that any
function $f \in L^p (D)$ can be approximated arbitrarily well in $L^p(D)$ norm
by a linear combination of the set of products $w_1 w_2$ of the elements
of the null-spaces $N_m$.

For example, if $L = \nabla^2$ then $N(\nabla^2)$ is the set of harmonic
functions, and the Laplacian has property $C$ if the set of products
$h_1 h_2$ of harmonic functions is total (complete) in
$L^p (D)$.

The notion of property $C$ has been introduced in \cite{R8}. It was developed
and widely used in \cite{R8} - \cite{R19}. It proved to be a very powerful
tool for a study of inverse problems \cite{R13} - \cite{R16}, \cite{R18} -
\cite{R19}.

Using property $C$ the author has proved in 1987 the uniqueness theorem for
3D inverse scattering problem with fixed-energy data \cite{R10}, \cite{R11},
\cite{R15}, uniqueness theorems for inverse problems of geophysics
\cite{R10}, \cite{R14}, \cite{R16}, and for many other inverse problems
\cite{R16}. The above problems have been open for several decades.

\subsection{} 
The third question that I will discuss, deals with approximation by entire
functions of exponential type. This question is quite simple but
the answer was not clear to engineers in the fifties,
it helped to understand the problem of
resolution ability of linear instruments \cite{R20}, \cite{R21}, and later
it turned to be useful in tomography \cite{R23}. This question in
applications is known as spectral extrapolation.

To formulate it, let us assume that $D \subset \R^n_x$ is a known
bounded domain, 
$$\widetilde f (\xi) := \int_D f(x) e^{i \xi \cdot x} dx := {\cal F} f,
\quad f(x) \in L^2 (D),  \eqno{(1.6)}$$
and assume that $\widetilde f (\xi)$ is known for $\xi \in \widetilde D$,
where $\widetilde D$ is a domain in $\R^n_\xi$.
{\it The question is: can one recover $f(x)$ from the knowledge of
$\widetilde f (\xi)$ in $\widetilde D$?}

Uniqueness of $f(x)$ with the data
$\{\widetilde f(\xi), \xi \in \widetilde D\}$
is immediate: $\widetilde f(\xi)$ is an entire function of exponential type
and if $\widetilde f(\xi) = 0$ in $\widetilde D$, then, by the analytic
continuation, $\widetilde f(\xi) \equiv 0$, and therefore $f(x)=0$.
Is it possible to derive an analytic formula for the recovery of $f(x)$
from $\{\widetilde f(\xi), \xi \in \widetilde D\}$?
It turns out that the answer is yes (\cite{R22} - \cite{R24}). Thus we give
an analytic formula for inversion of the Fourier transform
$\widetilde f(\xi)$ of a compactly supported function $f(x)$ from a
compact set $\widetilde D$.

From the point of view of approximation theory this problem is closely
related to the problem of approximation of a given function $h(\xi)$
by entire functions of exponential type whose Fourier transform has support
inside a given convex region. This region is
fixed but can be arbitrarily small.

In sections 2,3 and 4 the above three questions of approximation theory are
discussed in more detail, some of the results are formulated and some of them
are proved.

\section{Stable approximation of the derivative from noisy data.} 
In this section we formulate an answer to question 1.1. Denote
$\| f^{(1+a)} \| := m_{1+a},$ where $0 < a \leq 1$, and
$$\|f^{(1+a)} \| = \| f^\prime \|_\infty + \sup_{x,y \in \R}
  \frac{|f^\prime (x) - f^\prime (y)|}{|x-y|^a}.
  \eqno{(2.1)} $$

\begin{theorem} 
There does not exist an operator $T$ such that
$$\sup_{f \in {\cal K} (\delta, m_j)} \|Tf_\delta - f^\prime\|_\infty \leq
  \eta(\delta) \to 0 \hbox{\ as\ } \delta \to 0,
  \eqno{(2.2)}$$
if $j = 0$ or $j=1$. There exists such an operator if $j>1$. For example,
one can take $T=L_{\delta, j}$ where
$$L_{\delta, j} f_\delta :=
  \frac{f_\delta (x+h_j (\delta))-f_\delta (x-h_j (\delta))}{2h_j (\delta)},
  \quad h_j (\delta) := \left(\frac{\delta}{m_j(j-1)} \right)^{\frac{1}{j}},
  \eqno{(2.3)} $$
and then
$$\sup_{f \in {\cal K} (\delta, m_j)}
  \|L_{\delta, j} f_\delta - f^\prime \|_\infty \leq c_j
  \delta^{\frac{j-1}{j}}, \quad 1 < j \leq 2,
  \eqno{(2.4)} $$
$$c_j := \frac{j}{(j-1)^{\frac{j-1}{j}}} \quad m_j^{\frac{1}{j}}.
  \eqno{(2.5)}$$
\end{theorem}

\begin{proof} 
1. {\sf Nonexistence} of $T$ for $j=0$ and $j=1$.

Let $f_\delta (x) = 0$, $f_1 (x) := -\frac{mx(x-2h)}{2}, 0 \leq x \leq
2h$.
Extend
$f_1(x)$ on $\R$ so that $\| f_1^{(j)} \|_\infty =
  \sup_{0 \leq x \leq 2h} \| f_1^{(j)} \|, j = 0,1,2,$
and set $f_2 (x) := -f_1(x)$. Denote $(T f_\delta) (0) := b$.

One has
$$\| Tf_\delta - f_1^\prime \| \geq |(Tf_\delta) (0) - f^\prime_1 (0)| =
  |b-mh|, \eqno{(2.6)}$$
and
$$\| Tf_\delta - f^\prime_2 \| \geq |b+mh|. \eqno{(2.7)}$$

Thus, for $j=0$ and $j=1$, one has:
$$\gamma_j := \inf_T \sup_{f \in {\cal K} (\delta, m_j)}
  \| Tf_\delta - f^\prime \| \geq \inf_{b \in \R} \max
  \left[ |b-mh|, |b+mh| \right] = mh. \eqno{(2.8)}$$

Since $\|f_s - f_\delta\|_\infty \leq \delta, \quad s=1,2,$ and
$f_\delta = 0$, one gets
$$\| f_s \|_\infty = \frac{mh^2}{2} \leq \delta, \quad s=1,2.
  \eqno{(2.9)}$$
Take
$$h=\sqrt{\frac{2\delta}{m}}. \eqno{(2.10)}$$

Then
$m_0 = \| f_s \|_\infty = \delta, m_1 = \| f^\prime_s \|_\infty =
\sqrt{2\delta m}$,
so (2.8) yields
$$\gamma_0 = \sqrt{2\delta m} \to \infty \hbox{\ as\ } m \to \infty,
  \eqno{(2.11)}$$
and (2.2) does not hold if $j=0$.

If $j=1$, then (2.8) yields
$$\gamma_1 = \sqrt{2 \delta m} =m_1 > 0, \eqno{(2.12)}$$
and, again, (2.2) does not hold if $j=1$.

2. {\sf Existence} of $T$ for $j>1$.

If $j>1$, then the operator $T=L_{\delta, j}$, defined in (2.3) yields
estimate (2.4), so (2.2) holds with
$\eta(\delta) = c_j \delta^{\frac{j-1}{j}}$ and $c_j$ is defined in
(2.5). Indeed,
$$\| L_{\delta, j} f_\delta - f^\prime \|_\infty \leq \|L_{\delta, j}
  (f_\delta - f) \|_\infty + \|L_{\delta,j} f-f^\prime \|_\infty \leq
  \frac{\delta}{h} + m_jh^{j-1}. \eqno{(2.13)}$$
Minimizing the right-hand side of (2.13) with respect to $h>0$ for a fixed
$\delta >0$, one gets (2.4) and (2.5).

Theorem 2.1 is proved.
\end{proof}

\begin{remark}
If $j=2$, then $m_2=m$, where $m$ is the number introduced in the beginning
of the proof of Theorem 2.1, and using the Taylor formula one can get a
better estimate in the right-hand side of (2.13) for $j=2$, namely
$\|L_{\delta, 2} f_\delta - f^\prime \|_\infty \leq \frac{\delta}{h}
 + \frac{m_2h}{2}$.
Minimizing with respect to $h$, one gets
$h(\delta) = \sqrt{\frac{2\delta}{m_2}}$ and
$\min_{h>0} \left(\frac{\delta}{h} + \frac{m_2h}{2} \right) := \varepsilon
 (\delta) = \sqrt{2\delta m_2}$.
Now (2.8), with $m=m_2$, yields
$$\gamma_2 \geq \sqrt{2 \delta m_2} = \varepsilon (\delta).
  \eqno{(2.14)}$$
\end{remark}

Thus, we have obtained:
\begin{corollary} 
Among all linear and nonlinear operators $T$, the operator
$Tf= L_\delta f := \frac{f(x+h(\delta)) - f(x-h(\delta))}{2h(\delta)}$,
$h(\delta) = \sqrt{\frac{2\delta}{m_2}}$,
yields the best approximation of $f^\prime$, $f \in {\cal K}(\delta, m_2)$,
and
$$\gamma_2 := \inf_{T} \sup_{f \in {\cal K} (\delta, m_2)}
  \| Tf_\delta - f^\prime \|_\infty = \varepsilon (\delta) :=
  \sqrt{2\delta m_2}. \eqno{(2.15)}$$
\end{corollary}

\begin{proof}
We have proved that $\gamma_2 \geq \varepsilon (\delta)$.
If
$T=L_\delta$ then $\|L_\delta f_\delta - f^\prime \|_\infty \leq
 \varepsilon (\delta)$,
as follows from the Taylor's formula: if $h=\sqrt{\frac{2\delta}{m_2}}$,
then
$$\|L_\delta f_\delta -f^\prime \| \leq \frac{\delta}{h} +
  \frac{m_2h}{2} = \varepsilon (\delta), \eqno{(2.16)}$$
so $\gamma_2=\epsilon (\delta).$
\end{proof}

\section{Property $C$} 

\subsection{} 
In the introduction we have defined property $C$ for PDE. Is this property
generic or is it an exceptional one?

Let us show that this property is generic: a linear formal partial
differential operator with constant coefficients, in general, has
property $C$.In particular, the operators
$L=\nabla^2, L = i\partial_t - \nabla^2, L=\partial_t -\nabla^2$,
and $L=\partial^2_t - \nabla^2$, all have property $C$.

A necessary and sufficient condition for a pair $\{L_1, L_2\}$ of partial
differential operators to have property $C$ was found in \cite{R8} and
\cite{R26} (see also \cite{R16}).

Let us formulate this condition and use it to check that the four operators,
mentioned above, have property $C$.

Let $L_m u:= \sum_{|j| \leq J_m} a_{jm} D^j u(x)$, $m=1,2, a_{jm} = \const$,
$x \in \R^n$, $n\geq 2,$  $J_m \geq 1$.
Define the algebraic varieties
$${\cal L}_m := \{z : z \in \C^n, L_m (z) := 
\sum_{|j| \leq J_m} a_{jm} z^j = 0\},\quad
  m=1,2.$$

\begin{definition}
Let us write ${\cal L}_1 \nparallel {\cal L}_2$ if and only if there
exist at least one
point $z^{(1)} \in {\cal L}_1$ and at least one point
$z^{(2)} \in {\cal L}_2$, such
that the tangent spaces in $\C^n$ $T_m$ to ${\cal L}_m$ at the points
$z^{(m)}$, $m=1,2,$ are transversal, that is $T_1 \nparallel T_2$.
\end{definition}

\begin{remark3.2} 
A pair $\{{\cal L}_1, {\cal L}_2\}$ fails to have the property
${\cal L}_1 \nparallel {\cal L}_2$
if and only if ${\cal L}_1 \cup {\cal L}_2$ is a union of parallel
hyperplanes in $\C^n$.
\end{remark3.2}

\begin{theorem3.3} 
A pair $\{L_1, L_2\}$ of formal linear partial differential operators with
constant coefficients has property $C$ if and only if
${\cal L}_1 \nparallel {\cal L}_2$.
\end{theorem3.3}

\underline{Example 3.4}
Let
$L= \nabla^2$, then
${\cal L} = \{z:z \in \C^n, z_1^2 + z^2_2 + \dots z_n^2 = 0\}$.
It is clear that the tangent spaces to ${\cal L}$ at the points $(1,0,0)$ and
$(0,1,0)$ are transversal. So the Laplacian $L=\nabla^2$ does have property
$C$ (see Definition 1.2 in section 1). In other words, given an arbitrary
bounded domain $D \subset \R^n$ and an arbitrary function $f(x) \in L^p(D)$,
$p \geq 1$, for example, $p=2$, one can approximate $f(x)$ in the norm
of $L^P(D)$ by linear combinations of the product $h_1 h_2$ of harmonic
in $L^2(D)$ functions.

\underline{Example 3.5}
One can check similarly that the operators
$\partial_t - \nabla^2, i \partial_t - \nabla^2$ and
$\partial^2_t - \nabla^2$ have property $C$.

\nd{\bf Proof of Theorem 3.3}
We prove only the sufficiency and refer to \cite{R16} for the necessity.
Assume that
${\cal L}_1 \nparallel {\cal L}_2$. Note that
$e^{x \cdot z} \in N(L_m) := N_m$
if and only if $L_m (z) = 0, z \in \C^n$, that is $z \in {\cal L}_m$.
Suppose
$$
  \int_D f(x) w_1 w_2 dx = 0 \quad \forall w_m \in N_m,$$ 
then
$$
  F(z_1 + z_2) := 
  \int_D f(x) e^{x \cdot (z_1+z_2)} dx = 0 \quad
  \forall z_m \in {\cal L}_m.
  \eqno{(3.1)}$$
The function $F(z)$, defined in (3.1), is entire. It vanishes identically
if it vanishes on an open set in $\C^n$ (or $\R^n$). The set
$\{z_1 + z_2\}_{\forall z_m \in N_m}$ contains a ball in $\C^n$ if (and
only if) ${\cal L}_1 \nparallel {\cal L}_2$. Indeed, if $z_1$ runs though
the set ${\cal L}_1 \cap B(z^{(1)}, r)$, where
$B(z^{(1)}, r) := \{z: z\in \C^n, |z^{(1)} - z| < r\}$,
and $z_2$ runs through the set ${\cal L}_2 \cap B(z^{(2)}, r)$, then, for
all sufficiently small $r>0$, the set $\{z_1 + z_2\}$ contains a small
ball $B(\zeta, \rho)$, where $\zeta =z^{(1)} + z^{(2)}$, and $\rho >0$ is
a
sufficiently small number. To see this, note that $T_1$ has a basis
$h_1, \dots h_{n-1}$, which contains $n-1$ linearly independent vectors of
$\C^n$, and $T_2$ has a basis such that at least one of its vectors, call it
$h_n$, has a non-zero projection onto the normal to $T_1$, so that
$\{h_1, \dots, h_n\}$ are $n$ linearly independent vectors in $\C^n$. Their
linear combinations fill in a ball $B(\zeta, \rho)$. Since the vectors in
$T_m$ approximate well the vectors in
${\cal L}_m \cap B(z^{(m)}, r)$ if $r>0$ is sufficiently small, the set
$\{z_1 + z_2\}_{\forall z_m \in {\cal L}_m}$ contains a ball
$B(\zeta, \rho)$ if $\rho >0$ is sufficiently small. Therefore, condition (3.1)
implies $F(z) \equiv 0$, and so $f(x) = 0$. This means that the pair
$\{L_1, L_2\}$ has property $C$. We have proved:
$\{{\cal L}_1 \nparallel {\cal L}_2\} \Rightarrow \left\{ \{L_1, L_2\}
  \hbox{\ has\ } \hbox{property\ } C \right\}$.
\qed

How does one prove property $C$ for a pair
$\{L_1, L_2\} := \{\nabla^2 + k^2 - q_1 (x), \nabla^2 + k^2 - q_2 (x)\}$
of the Schr\"odinger operators, where $k = \const \geq 0$ and 
$q_m (x) \in L^2_{loc} (\R^n)$ are
some real-valued, compactly supported functions?

One way to do it \cite{R16} is to use the existence of the elements
$\psi_m \in N_m := \{w:[\nabla^2 + k^2 - q_m (x)] w=0 \hbox{\ in\ } \R^n\}$
which are of the form
$$\psi_m (x, \theta) = e^{ik \theta \cdot x} [1 + R_m(x, \theta)], \quad
  k > 0, \eqno{(3.2)}$$
where $\theta \in M := \{ \theta : \theta \in \C^n, \theta \cdot \theta =1\}.$
Here $\theta \cdot w := \sum^n_{j=1} \theta_j w_j$, (note that there is no
complex
conjugation above $w_j$), the variety $M$ is noncompact, and \cite{R16}
$$\|R_m (x, \theta) \|_{L^\infty (D)} \leq c
  \left(\frac{\ln |\theta|}{|\theta|} \right)^{\frac{1}{2}}, \quad
  \theta \in M, \quad
  |\theta| \to \infty, \quad  m=1,2, \eqno{(3.3)}$$
where $c = \const >0$ does not depend on $\theta$, $c$ depends on $D$
and on $\| q \|_{L^\infty(B_a)}$, where $q = \overline q, q=0$
for $|x| >a$, and $D \subset \R^n$ is an arbitrary bounded domain. Also
$$\| R_m (x, \theta) \|_{L^2(D)} \leq \frac{c}{|\theta|}, \quad
        \theta \in M, \quad |\theta| \to \infty, \quad m=1,2.
        \eqno{(3.4)}$$
It is easy to check that for any $\xi \in \R^n, n \geq 3$,
and any $k>0,$ one can find
(many)
$\theta_1$ and $\theta_2$ such that
$$k(\theta_1 + \theta_2) = \xi,\quad |\theta_1| \to \infty,
 \quad  \theta_1, \theta_2 \in M, \quad n \geq 3. \eqno{(3.5)}$$
Therefore, using (3.5) and (3.3), one gets:
$$\lim_{\substack{|\theta_1| \to \infty \\
                \theta_1 + \theta_2 = \xi,\,\, \theta_1 ,  \theta_2 \in
M}}
        \psi_1 \psi_2 = e^{i \xi \cdot x},
        \eqno{(3.6)}$$
Since the set $\{e^{i \xi \cdot x}\}_{\forall \xi \in \R^n}$
is total in $L^p(D)$, it follows that the pair $\{L_1, L_2\}$
of the Schr\"odinger operators under the above assumptions does have
property $C$.

\subsection{} 
Consider the following problem of approximation theory \cite{R19}.

Let $k=1$ (without loss of generality), $\alpha \in S^2$ (the unit sphere
in $\R^3$), and $u := u(x, \alpha)$ be the scattering solution that is, an
element of $N(\nabla^2 + 1- q(x))$ which solves the problem:
$$\left[\nabla^2 + 1-q(x) \right]u= 0 \hbox{\ in\ } \R^3, \eqno{(3.5)}$$
$$u= \exp (i \alpha \cdot x) + A(\alpha^\prime, \alpha)
  \frac{e^{i|x|}}{|x|} + o \left(\frac{1}{|x|} \right),\quad |x| \to
\infty,\quad
  \alpha^\prime := \frac{x}{|x|}, \eqno{(3.6)}$$
where $\alpha \in S^2$ is given.

Let $w \in N(\nabla^2 + 1-q(x)) := N(L)$ be arbitrary, $w \in H^2_{loc}$,
$H^l$ is the Sobolev space. The problem is: is it possible to approximate
$w$ in $L^2(D)$ with an arbitrary accuracy by a linear combination of the
scattering solutions $u(x, \alpha)$? In other words, given an arbitrary
small number $\varepsilon > 0$ and an arbitrary fixed, bounded, 
homeomorphic
to a ball, Lipschitz domain $D \subset \R^n$, can one find
$v_\varepsilon (\alpha) \in L^2 (S^2)$ such that
$$\| w - \int_{S^2} u(x, \alpha) v_\varepsilon (\alpha) d \alpha \|_{L^2(D)}
  \leq \varepsilon ? \eqno{(3.7)}$$
If yes, what is the behavior of $\|v_\varepsilon (\alpha) \|_{L^2(S^2)}$
as $\varepsilon \to 0$, if $w = \psi(x,\theta)$ where $\psi$ is the special
solution (3.2) - (3.3), $\theta \in M, Im \theta \neq 0$?

The answer to the first question is yes. A proof \cite{R16} can go as follows.
If (3.7) is false, then one may assume that $w \in N(L)$ is such that
$$\int_D \overline w \left(\int_{S^2} u(x, \alpha) v (\alpha) d \alpha
  \right) dx = 0 \quad \forall v \in L^2(S^2). \eqno{(3.8)}$$
This implies
$$\int_D \overline w u(x, \alpha) dx = 0 \quad \forall \alpha \in S^2.
  \eqno{(3.9)}$$
From (3.9) and formula (5) on p. 46 in \cite{R16} one concludes:
$$v(y) := \int_D \overline w G(x,y) dx = 0 \quad \forall y \in D^\prime :=
  \R^3 \backslash D, \eqno{(3.10)}$$
where $G(x,y)$ is the Green function of $L$:
$$\left[\nabla^2 + 1 - q(x) \right] G(x,y) = -\delta(x,y) \hbox{\ in\ }
  \R^3, \eqno{(3.11)}$$
$$\lim_{r \to \infty} \int_{|x| = r} \left| \frac{\partial G}{\partial |x|}
  - iG \right|^2 ds = 0. \eqno{(3.12)}$$

Note, that in [\cite{R16}, p. 46] formula (5) is:
$$G(x,y, k) = \frac{e^{i|y|}}{4\pi |y|} u (x, \alpha) +
  o \left(\frac{1}{|y|} \right),\quad |y| \to \infty, \quad \frac{y}{|y|}
= -\alpha.
  \eqno{(3.13)}$$
From (3.10) and (3.11) it follows that
$$\left[\nabla^2 + 1 - q(x) \right] v(x) = -\overline w (x) \hbox{\ in\ } D,
  \eqno{(3.14)} $$
$$v = v_N = 0 \hbox{\ on\ } S := \partial D, \eqno{(3.15)}$$
where $N$ is the outer unit normal to $S$, and (3.15) holds because
$v = 0$ in $D^\prime$ and $v \in H^2_{loc}$ by elliptic regularity if
$q \in L^2_{loc}$. Multiply (3.14) by $w$, integrate over $D$ and then, on
the left, by parts, using (3.15), and get:
$$\int_D |w|^2 dx =0. \eqno{(3.16)}$$

Thus, $w=0$ and (3.7) is proved.
\qed

Let us prove that $\| v_\varepsilon (\alpha) \|_{L^2(S^2)} \to \infty$
as $\varepsilon \to 0$. Assuming
$\|v_\varepsilon (\alpha) \| \leq c \quad \forall \varepsilon
 \in (0, \varepsilon_0)$, where $\varepsilon_0 > 0$ is some number, one can
select a weakly convergent in $L^2(S^2)$ sequence
$v_n(\alpha) \to v(\alpha), n \to \infty$, and pass to the limit in (3.7)
with $w = \psi(x, \theta)$, to get
$$\| \psi (x, \theta) - \int_{S^2} u(x, \alpha) \nu(\alpha)
   d \alpha \|_{L^2(D)} = 0. \eqno{(3.17)}$$
The function
$U(x) := \int_{S^2} u(x, \alpha) \nu (\alpha) d \alpha \in N(L)$ and
$\|U(x)\|_{L^\infty(\R^n)} < \infty$
(since $\sup_{x \in \R^n, \alpha \in S^2}$
$|u(x, \alpha)| < \infty$).
Since (3.17)
implies that $\psi (x, \theta) = U(x)$ in $D$, and both $\psi(x, \theta)$
and $U(x)$ solve equation (3.5), the unique continuation principle for the
solutions of the elliptic equation (3.5) implies
$U(x) = \psi(x, \theta)$ in $\R^3$. This is a contradiction since
$\psi(x, \theta)$ grows exponentially as $|x| \to \infty$ in certain
directions because $Im \theta \neq 0$ (see formula (3.2)).

We have proved that if $w = \psi(x, \theta), Im \theta \neq 0$, then
$\| \nu_\varepsilon (\alpha) \|_{L^2(S^2)} \to \infty$ as
$\varepsilon \to 0$, where $\nu_\varepsilon (\alpha)$ is the function from
(3.7).

For example, if $q(x) = 0$ and $k=1$ then
$\psi(x, \theta) = e^{i \theta \cdot x}$ and
$u(x, \alpha) = e^{i\alpha \cdot x}$. So, if $Im \theta \neq 0, \theta \in M$,
and $\|e^{i\theta \cdot x} -\int_{S^2} \nu_\varepsilon (\alpha)
     e^{i\alpha \cdot} \|_{L^2(D)} < \varepsilon$,
then $\| \nu_\varepsilon (\alpha) \|_{L^2(S^2)} \to \infty$ as
$\varepsilon \to 0$. It is interesting to estimate the rate of growth of
$$\inf_{\substack{\nu \in L^2(S^2) \\
                \|e^{i\theta \cdot x} - \int_{S^2} \nu(\alpha)
                e^{i\alpha \cdot x} d \alpha \|_{L^2(D)} < \varepsilon}}
   \| \nu (\alpha) \|_{L^2(S^2)}.$$
This is done in \cite{R19} (and \cite{R27}).

Property $C$ for ordinary differential equations is defined, proved and
applied to many inverse problems in \cite{R17}.

\section{Approximation by entire functions of exponential type.} 
Let $f(x) \in L^2(B_a)$, $B_a := \{x ; x \in \R^n_x, |x| \leq a\}$,
$a >0$ is a fixed number, $f(x)=0$ for $|x|>a$, and
$$\widetilde f(\xi) = \int_{B_a} f(x) e^{i \xi \cdot x} dx. \eqno{(4.1)}$$
Assume that $\widetilde f(\xi)$ is known for all
$\xi \in \widetilde D \subset \R^n_\xi$, where $\widetilde D$ is a
(bounded) domain.

The problem is to find $\widetilde f (\xi)$ for all $\xi \in \R^n_\xi$,
(this is called spectral extrapolation), or, equivalently, to find
$f(x)$ (this is called inversion of the Fourier transform $\widetilde f(\xi)$
of a compactly supported function $f(x)$, $\supp f(x) \subset B_a$, from
a compact $\widetilde D$).

In applications the above problem is also of interest in the case when
$\widetilde D$ is not necessarily bounded. For example, in tomography
$\widetilde D$ may be a union of two infinite cones (the limited-angle data).

In the fifties and sixties there was an extensive discussion in the
literature concerning the resolution ability of linear instruments. According to
the theory of the formation of optical images, the image of a bright point,
which one obtains when the light, issued by this point, is diffracted on a
circular hole in a plane screen, is the Fourier transform
$\widetilde f(\xi)$ of the function $f(x)$ describing the light
distribution on the
circular hole $B_a$, the two-dimensional ball. This Fourier transform is an
entire function of exponential type. One says that the resolution ability
of a linear instrument (system) can be increased without a limit if the
Fourier transform of the function $f(x)$, describing the light distribution
on $B_a,$ can approximate the delta-function $\delta(\xi)$ with an
arbitrary accuracy. The above definition is not very precise, but it is
usual in applications and can be made precise: it is sufficient to specify
the metric in which the delta-function is approximated. For our purposes,
let us take a delta-type sequence $\delta_j (\xi)$ of continuous functions
which is defined by the requirements:
$\left| \int_{\widetilde D} \delta_j (\xi) d \xi \right| \leq c$,
where the constant $c$ does not depend on $\widetilde D$ and $j$, and
$$\lim_{j \to \infty} \int_{\widetilde D} \delta_j (\xi) d \xi =
 \begin{cases} 1 \hbox{\ if\ } 0 \in \widetilde D, \\
               0 \hbox{\ if\ } 0 \not\in \widetilde D.
              \end{cases}$$

{\it The approximation problem is: can one approximate an arbitrary continuous (or
$L^1(\widetilde D)$) function $g(\xi)$ in an arbitrary fixed bounded
domain $\widetilde D \subset \R^n_\xi$ by an entire function of
exponential type
$\widetilde f (\xi) = \int_{B_a} f(x) e^{i \xi \cdot x} dx$, where
$a> 0$ is an
arbitrary small number?}

The engineers discussed this question in a different form: can one
transmit with an arbitrary accuracy a high-frequency signal
$g(\xi)$ by using low-frequency signals $\widetilde f(\xi)$?
The smallness of $a$ means that the "spectrum" $f(x)$
of the signal $\widetilde f(\xi)$ contains only
"low spatial frequencies".

From the mathematical point of view the answer is nearly obvious: yes.
The proof is very simple: if an approximation with an arbitrary accuracy
were impossible, then
$$0 = \int_{\widetilde D} g(\xi) \left( \int_{B_a} e^{i \xi \cdot x}
  f(x) dx \right) d \xi \quad \forall f \in L^2(B_a).$$

This implies the relation
$$0 = \int_{\widetilde D} g(\xi) e^{i \xi \cdot x} d \xi \quad
  \forall x \in B_a.$$
Since $\widetilde D$ is a bounded domain, the integral above is an entire
function of $x \in \C^n$ which vanishes in a ball $B_a$. Therefore this
function vanishes identically and consequently $g(\xi) \equiv 0$. This
contradiction
proves that the approximation of an arbitrary
$g(\xi) \in L^2 (\widetilde D)$ by the entire functions
$\widetilde f (\xi) = \int_{B_a} f(x) e^{i \xi \cdot x} dx$
is possible with an arbitrary accuracy in $L^2(\widetilde D)$.

Now let us turn to another question: how does one derive an analytic formula
for finding $f(x)$ if $\widetilde f (\xi)$ is given in
$\widetilde D$?

In other words, how does one invert analytically the Fourier transform
$\widetilde f(\xi)$ of a compactly supported function $f(x)$,
$\supp f \subset B_a$, from a compact $\widetilde D$?

We discuss this question below, but first let us discuss the notion of
apodization, which was a hot topic at the end of the sixties. Apodization
is a method to increase the resolution ability of a linear optical system
(instrument) by putting a suitable mask on the outer pupil of the
instrument. Mathematically one deals with an approximation problem: by
choosing a mask $g(x)$, which transforms the function $f(x)$ on the outer
pupil of the instrument into a function $g(x) f(x)$, one wishes to
change the image $\widetilde f(\xi)$ on the image plane to an image
$\delta_j (\xi)$ which is close to the delta-function $\delta(\xi)$,
and therefore increase the resolution ability of the instrument. That the
resolution ability can be increased without a limit (only in the absence
of noise!) follows from the above argument: one can choose $g(x)$ so that
$\widetilde{g(x) f(x)}$ will approximate arbitrarily accurately
$\delta_j(\xi)$, which, in turn, approximates $\delta (\xi)$ arbitrarily
accurately.

This conclusion contradicts to the usual intuitive idea according to
which one cannot resolve details smaller than the wavelength.

In fact, if there is no noise, one can, in principle, increase resolution
ability without a limit (superdirectivity in the antenna theory), but
since the noise is always present, in practice there is a limit to the
possible increase of the resolution ability.

Let us turn to the analytic formula for the approximation by entire
functions and for the inversion of the Fourier transform of a compactly
supported function from a compact $\widetilde D$.

Multiply (4.1) by $(2 \pi)^{-n} \widetilde{\delta_j} (\xi) e^{-i \xi \cdot x}$
and integrate over $\widetilde D$ to get
$$f_j (x) = \int_{B_a} f(y) \delta_j (x-y) dy = \frac{1}{(2 \pi)^n}
  \int_{\widetilde D} \widetilde f (\xi) \widetilde{\delta_j}
  e^{-i \xi \cdot x} d \xi \eqno{(4.2)}$$
where $\widetilde{\delta_j} (\xi)$ is the Fourier transform of
$\delta_j (x)$.

Let us choose $\delta_j (x)$ so that it will be a delta-type sequence (in the
sense defined above). In this case $f_j(x)$ approximates $f(x)$
arbitrarily accurately:
$$\| f - f_j (x) \| \to 0 \hbox{\ as\ } j \to \infty, \eqno{(4.3)}$$
where the norm $\| \cdot \|$ is  $L^2(B_a)$ norm if $f \in L^2(B_a)$,
and $C(B_a)$-norm if $f \in C(B_a)$.

If $\| f \|_{C^1 (B_a)} \leq m_1$, then
$$\| f-f_j (x) \|_{C(B_a)} \leq c m_1 j^{-\frac{1}{2}}. \eqno{(4.4)}$$

The conclusions (4.3) and (4.4) hold if, for example, 
$$\delta_j (x) := P_j (|x|^2) \left({\cal F}^{-1} \widetilde h \right)
  (x), \eqno{(4.5)}$$
where
$$P_j (r) := \left(\frac{j}{4 \pi a_1^2} \right)^{\frac{n}{2}}
  \left(1-\frac{r}{4a^2_1} \right)^j, \quad 0 \leq r \leq a,\quad
  a_1 > a, \eqno{(4.6)}$$
$$\widetilde h (\xi) \in C^\infty_0 (\widetilde D), \quad \frac{1}{(2\pi)^n}
  \int_{\widetilde D} \widetilde h (\xi) d \xi = 1. \eqno{(4.7)}$$

\begin{theorem} 
If (4.5)-(4.7) hold, then the sequence $f_j(x)$, defined in (4.2),
satisfies (4.3) and (4.4).
\end{theorem}

Thus, formula (4.2):
$$f(x) = {\cal F}^{-1} \left[ \widetilde f (\xi) \widetilde{\delta_j}
  (\xi) \right], \eqno{(4.8)} $$
where $\widetilde{\delta_j} (\xi) := {\cal F} \delta_j (x)$, and
$\delta_j (x)$ is defined by formulas (4.5)-(4.7), is an inversion formula
for the Fourier transform $\widetilde f(\xi)$ of a compactly supported
function $f(x)$ from a compact $\widetilde D$ in the sense (4.3). A proof
of a Theorem similar to 4.1 has been originally published in \cite{R22}.

In \cite{R20}, \cite{R21} apodization theory and resolution ability are
discussed. In \cite{R24} a one-dimensional analog of Theorem 4.1 is given.
In this analog one can choose analytically explicitly a function similar
to $\widetilde h(\xi)$.

Proof of Theorem 4.1 is given in [\cite{R23}, pp. 260-263].

Let us explain a possible application of Theorem 4.1 to the limited-angle
data in tomography. The problem is: let
$\widehat f (\alpha, p) := \int_{l_{\alpha p}} f(x) ds$,
where $l_{\alpha, p}$ is the Radon transform of $f(x)$. Given
$\widehat f (\alpha, p)$ for all $p \in \R$ and all $\alpha \in K$, where
$K$ is an open proper subset of $S^2$, find $f(x)$, assuming
$\supp f \subset B_a$. It is well known \cite{R23}, that
$$\int^\infty_{-\infty} \widehat f (\alpha, p) e^{i pt} dp =
  \widetilde f (t \alpha), \quad t \in \R, \quad t \alpha := \xi.$$
Therefore, if one knows $\widehat f (\alpha, p)$ for all $\alpha \in K$
and all $p \in \R$, then one knows $\widetilde f (\xi)$ for all $\xi$ in a
cone $K \times \R$. Now Theorem 4.1 is applicable for finding $f(x)$
given $\widehat f (\alpha, p)$ for $\alpha \in K$ and $p \in \R$.


\begin{thebibliography}{10}

\bibitem{DI}
Dolgopolova T., Ivanov V.,
On numberical differentiation, Jour. of Comput. Math. and Math. Phys.,
6, N3, (1966), 570-576.

\bibitem{P}
Phillips D.L.,
A technique for numerical solution of certain integral equations of the first
kind, J. Assoc. Comput. Math., 9, N1, (1962), 84-97.

\bibitem{R1}
Ramm, A.G.,
On numerical differentiation. Mathem., Izvestija
vuzov, 11, 1968, 131-135. Math. Rev. 40 \#5130.

\bibitem{R2}
Ramm, A.G.,
{\bf Scattering by obstacles, D.Reidel, Dordrecht, 1986, pp.1-442.}

\bibitem{R3}
Ramm, A.G.,
{\bf Random fields estimation theory, Longman
Scientific and Wiley, New York, 1990.}

\bibitem{R4}
Ramm, A.G.,
Stable solutions of some ill-posed problems, Math. Meth. in the
appl. Sci. 3, (1981), 336-363.

\bibitem{R5}
Ramm, A.G.,
Simplified optimal differentiators. Radiotech.i
Electron.17, (1972), 1325-1328; English translation pp.1034-1037.

\bibitem{R6}
Ramm, A.G.,
Inequalities for the derivatives, Math. Ineq. and
Appl., 3, N1, (2000), 129-132.

\bibitem{R7}
Ramm, A.G., A. B. Smirnova,
On stable numerical differentiation, Mathem. of Computation,
70,(2001), 1131-1153.


\bibitem{R8}
Ramm, A.G.,
On completeness of the products of harmonic
functions, Proc. A.M.S., 99, (1986), 253-256.

\bibitem{R9}
Ramm, A.G.,
Inverse scattering for geophysical problems. Phys.
Letters. 99A, (1983), 258-260.

\bibitem{R10}
Ramm, A.G.,
Completeness of the products of solutions to PDE
and uniqueness theorems in inverse scattering,
Inverse problems, 3, (1987), L77-L82


\bibitem{R11}
Ramm, A.G.,
Recovery of the potential from fixed energy
scattering data. Inverse Problems, 4, (1988),
877-886; 5, (1989) 255.

\bibitem{R12}
Ramm, A.G.,
Multidimensional inverse problems and completeness
of the products of solutions to PDE. J. Math. Anal.
Appl. 134, 1, (1988), 211-253; 139, (1989) 302.

\bibitem{R13}
Ramm, A.G.,
Multidimensional inverse scattering problems and
completeness of the products of solutions to
homogeneous PDE. Zeitschr. f. angew. Math. u.
Mech., 69, (1989) N4, T13-T22.

\bibitem{R14}
Ramm, A.G.,
Completeness of the products of solutions of PDE
and inverse problems, Inverse Probl.6,
(1990), 643-664.

\bibitem{R15}
Ramm, A.G.,
Uniqueness theorems for multidimensional inverse
problems with unbounded coefficients.
J. Math. Anal. Appl. 136, (1988), 568-574.

\bibitem{R16}
Ramm, A.G.,
{\bf  Multidimensional inverse scattering
problems, Longman/Wiley, New York, 1992, pp.1-385.}

\bibitem{R17}
Ramm, A.G.,
Property C for ODE and applications to inverse problems,
in the book "Operator Theory and Its Applications",
Amer. Math. Soc., Fields Institute Communications vol. 25,
(2000), pp.15-75,
Providence, RI.
(editors A.G.Ramm, P.N.Shivakumar, A.V.Strauss).

\bibitem{R18}
Ramm, A.G.,
Stability of solutions to inverse scattering problems
with fixed-energy data, Rend. del seminario Mat. e Fisico di Milano,
(2001), pp.135-211.

\bibitem{R19}
Ramm, A.G.,
Stability estimates in inverse scattering, Acta
 Appl. Math., 28, N1, (1992), 1-42.

\bibitem{R20}
Ramm, A.G.,
Apodization theory II. Opt. and Spectroscopy, 29,
(1970), 390-394.

\bibitem{R21}
Ramm, A.G.,
Increasing of the resolution ability of the optical
instruments by means of apodization,  Opt. and
Spectroscopy, 29,
(1970), 594-599.

\bibitem{R22}
Ramm, A.G.,
Approximation by entire functions, Mathematics,
Izv. vusov, 10, (1978), 72-76.

\bibitem{R23}
Ramm, A.G.,  A. I. Katsevich,
{\bf The Radon Transform and Local Tomography, CRC
Press, Boca Raton 1996, pp.1-503.} 

\bibitem{R24}
Ramm, A.G.,
Signal estimation from incomplete data, J. Math.
 Anal. Appl., 125 (1987), 267-271.


\bibitem{R25}
Ramm, A.G.,
On simultaneous approximation of a function and its
derivative by interpolation polynomials. Bull.
Lond. Math. Soc. 9, 1977, 283-288.

\bibitem{R26}
Ramm, A.G.,
Necessary and sufficient condition for a PDE to
have property C, J. Math. Anal. Appl.156, (1991), 505-509.

\bibitem{R27}
Ramm, A.G.,
{\bf  Multidimensional inverse scattering problems, Mir
Publishers, Moscow, 1994, pp.1-496. (Russian translation of
the expanded monograph [18]).}


\bibitem{T1}
Tikhonov, A.N.,
Solving ill-posed problems and regularization,
Donlady Acad. Sct. USSR, 157, N3, (1963), 501-504.


\end{thebibliography}
\end{document}